\documentclass[11pt]{article}

\usepackage[a4paper,margin=1in]{geometry}
\usepackage{amsmath,amssymb,amsthm,mathtools}
\usepackage{enumitem}
\usepackage{microtype}
\usepackage{booktabs}
\usepackage[hidelinks]{hyperref}

\newtheorem{maintheorem}{Theorem}

\newtheorem{maincorollary}[maintheorem]{Corollary}
\newtheorem{theorem}{Theorem}[section]

\newtheorem{lemma}[theorem]{Lemma}

\theoremstyle{definition}
\newtheorem{definition}[theorem]{Definition}
\theoremstyle{remark}
\newtheorem{remark}[theorem]{Remark}

\numberwithin{equation}{section}

\newcommand{\R}{\mathbb{R}}
\newcommand{\N}{\mathbb{N}}
\newcommand{\Tmax}{T_{\max}}
\newcommand{\Ts}{T^{*}}
\newcommand{\dd}{\,\mathrm{d}}
\newcommand{\eps}{\varepsilon}

\title{Instantaneous complete blow-up for the doubly critical\\
semilinear heat equation}
\author{Berikbol T. Torebek\\[2mm]
\small Institute of Mathematics and Mathematical Modeling, Almaty, Kazakhstan\\
\small SDU University, Kaskelen, Kazakhstan\\
\small \texttt{torebek@math.kz}}
\date{}

\begin{document}
\maketitle

\begin{abstract}
We study $u_t-\Delta u=u^{1+2/d}$ in a smooth bounded domain with
Dirichlet boundary condition and nonnegative data $u_0\in L^1$. Let
$\Ts(u_0)$ be the maximal existence time of nonnegative integral
solutions. We show that the solutions of the problems with truncated
nonlinearity $\min\{u^p,n\}$ converge to the minimal integral solution up
to $\Ts(u_0)$ and diverge everywhere afterwards, and that
$\Ts(\kappa u_0)\le\lim_{n\to\infty}\Tmax(\min\{u_0,n\})\le\Ts(u_0)$ for
every $\kappa>1$. Hence nonexistence of a local integral solution forces
instantaneous complete blow-up, and the datum constructed by Brezis and
Cazenave [J.
Anal. Math. 68 (1996), 277–304] gives a positive answer to their Open Problem~5. We also prove
the critical concentration estimate
$\Ts(u_0)\le C_d\rho^2\exp(C_d\mu^{-2/d})$ for data carrying mass $\mu$ in
a ball of radius $\rho$, with optimal exponent. Consequences include
optimal collapse rates for truncations, a dichotomy at $L\log^{d/2}L$, and
ill-posedness at every nonnegative datum.
\end{abstract}



\section{Introduction}

\subsection{The problem}

Let $\Omega\subset\R^d$, $d\ge1$, be a bounded domain with smooth
boundary, and consider
\begin{equation}\label{P}
\begin{cases}
 u_t-\Delta u=u^p & \text{in }(0,T)\times\Omega,\\
 u=0 & \text{on }(0,T)\times\partial\Omega,\\
 u(0)=u_0\ge0 & \text{in }\Omega,
\end{cases}
\qquad p=1+\frac2d .
\end{equation}
For $q_c:=d(p-1)/2$ the space $L^{q_c}$ is invariant under the scaling
$u_\lambda(t,x)=\lambda^{2/(p-1)}u(\lambda^2t,\lambda x)$, and
$p=1+2/d$ is exactly the case $q_c=1$, in which this scale-invariant
space is $L^1$. Following Brezis and Cazenave \cite[Section~7.5]{BC96}
we call this case \emph{doubly critical}.\footnote{The same name is also
used for the case $q=q_c=p$, that is $p=d/(d-2)$, in which uniqueness in
$C([0,T];L^{q_c})$ fails; see e.g.\ \cite{MT03}.}

For $q>q_c$, $q\ge1$, or $q=q_c>1$, problem~\eqref{P} has, for every
$u_0\in L^q(\Omega)$, a local solution in
$C([0,T];L^q(\Omega))\cap L^\infty_{\rm loc}((0,T];L^\infty(\Omega))$,
which is classical for $t>0$ and unique in this class \cite{We80,BC96}.
In the doubly critical case Brezis and Cazenave proved
\cite[Theorem~11]{BC96} that there is $0\le u_0\in L^1(\Omega)$ for
which \eqref{P} has no nonnegative solution in
$C([0,T];L^1(\Omega))\cap L^\infty_{\rm loc}((0,T);L^\infty(\Omega))$ for
any $T>0$, and they formulated a list of open problems
\cite[Section~7.5]{BC96}. Open Problem~1 was answered by Arrieta and
Rodr\'{\i}guez-Bernal \cite{AR04}, Open Problem~2 by Celik and Zhou
\cite{CZ03} (see also \cite{LRSV16}), and Open Problems~3 and 4, on
nonuniqueness in $C([0,T];L^1)$ and on failure of positivity
preservation, in \cite{To}. This paper is concerned with Open
Problem~5.

For $0\le w_0\in L^\infty(\Omega)$, problem \eqref{P} has a unique
maximal classical solution, with maximal existence time
$\Tmax(w_0)\in(0,\infty]$. For $n\in\N$ put
\[
 u_{0,n}:=\min\{u_0,n\},\qquad g_n(s):=\min\{s^p,n\}\quad(s\ge0),
\]
and let $U_n$ be the global solution of the problem with truncated
nonlinearity
\begin{equation}\label{Pn}
\begin{cases}
 \partial_tU_n-\Delta U_n=g_n(U_n) & \text{in }(0,\infty)\times\Omega,\\
 U_n=0 & \text{on }(0,\infty)\times\partial\Omega,\\
 U_n(0)=u_0 & \text{in }\Omega
\end{cases}
\end{equation}
(Lemma~\ref{lem:Un}).

\medskip
\noindent\textbf{Open Problem 5} \cite[Section~7.5]{BC96}.
\emph{Is there $0\le u_0\in L^1(\Omega)$ such that
$\Tmax(u_{0,n})\to0$ as $n\to\infty$? Alternatively, is there
$0\le u_0\in L^1(\Omega)$ such that $U_n(t,x)\to+\infty$ as $n\to\infty$
for all $t>0$ and all $x\in\Omega$?}

\medskip
Brezis and Cazenave pose their problems for $u_t-\Delta u=|u|^{p-1}u$,
with the odd truncation $\min\{|t|^p,n\}\,\mathrm{sign}\,t$; for
$u_0\ge0$ the truncated solutions are nonnegative, so this coincides
with the formulation above. Their Theorem~11 was presented as evidence
for a positive answer. The difficulty in deducing it from nonexistence is
that the truncated solutions might exist on a common time interval and
diverge only in the limit $n\to\infty$; we show that this cannot happen.

To formulate the results, let $S(t)$ be the Dirichlet heat semigroup
(Section~\ref{sec:prelim}).

\begin{definition}\label{def:integral}
Let $0\le u_0\in L^1(\Omega)$ and $T\in(0,\infty]$. A measurable
function $u:\Omega\times(0,T)\to[0,\infty]$ is a (nonnegative)
\emph{integral solution} of \eqref{P} on $(0,T)$ if $u<\infty$ almost
everywhere in $\Omega\times(0,T)$ and
\begin{equation}\label{integral}
 u(t)=S(t)u_0+\int_0^tS(t-s)u(s)^p\dd s
 \qquad\text{a.e. in }\Omega\times(0,T).
\end{equation}
We set
\[
 \Ts(u_0):=\sup\big\{T>0:\ \eqref{P}\text{ has a nonnegative integral
 solution on }(0,T)\big\}\in[0,\infty],
\]
with $\sup\emptyset:=0$.
\end{definition}

This is the notion of \cite[Definition~3.1]{LRSV16}, the weakest one for
which \eqref{P} makes sense. Classical solutions with bounded data and
nonnegative solutions in the class of \cite[Theorem~11]{BC96} are
integral solutions (Remark~\ref{rem:classes}); in particular
$\Tmax(w_0)\le\Ts(w_0)$ for bounded $w_0$. Open Problem~5 concerns the
following three properties of a datum $0\le u_0\in L^1(\Omega)$:
\begin{enumerate}[label=\textnormal{(\roman*)}]
\item $\Tmax(\min\{u_0,n\})\to0$ as $n\to\infty$;
\item $U_n(t,x)\to+\infty$ as $n\to\infty$, for every $t>0$ and every
 $x\in\Omega$;
\item $\Ts(u_0)=0$, that is, for no $T>0$ does \eqref{P} have a
 nonnegative integral solution on $(0,T)$.
\end{enumerate}
Finally, we say that \emph{bounded solutions blow up completely} if
\begin{equation}\label{CB}\tag{CB}
 \Ts(w_0)=\Tmax(w_0)\qquad\text{for every }0\le w_0\in L^\infty(\Omega),
\end{equation}
that is, bounded solutions cannot be continued as integral solutions
beyond $\Tmax$. This is equivalent to complete blow-up in the sense of
Baras and Cohen, and it holds if $\Omega$ is convex
(Remark~\ref{rem:BCo}).

\subsection{Main results}

\begin{maintheorem}[Answer to Open Problem~5]\label{cor:OP5}
Let $0\le u_0\in L^1(\Omega)$.
\begin{enumerate}[label=\textnormal{(\alph*)}]
\item \textnormal{(ii)}$\Leftrightarrow$\textnormal{(iii)}$\Rightarrow$\textnormal{(i)}.
\item If \textnormal{(i)} holds, then $\Ts(\kappa u_0)=0$ for every
 $\kappa>1$.
\item If \eqref{CB} holds, in particular if $\Omega$ is convex, then
 \textnormal{(i)--(iii)} are equivalent.
\item The datum constructed in the proof of \cite[Theorem~11]{BC96}
 satisfies \textnormal{(i)--(iii)}.
\end{enumerate}
In particular, in every smooth bounded domain both formulations of Open
Problem~5 have a positive answer.
\end{maintheorem}

Parts (a), (b) and (d) hold in every smooth bounded domain. Theorem~A is
a consequence of the following description of the truncated problems.

\begin{maintheorem}[Minimal solution and truncations]\label{thm:main}
Let $0\le u_0\in L^1(\Omega)$ and $U:=\lim_{n\to\infty}U_n$
(Lemma~\ref{lem:U}).
\begin{enumerate}[label=\textnormal{(\alph*)}]
\item $U(t,x)=+\infty$ for all $t>\Ts(u_0)$ and all $x\in\Omega$. If
 $\Ts(u_0)>0$, then $U$ is a nonnegative integral solution of \eqref{P}
 on $(0,\Ts(u_0))$, $U(t)\in L^1_{\rm loc}(\Omega)$ for
 $0<t<\Ts(u_0)$, and $U\le v$ a.e. for every nonnegative integral
 solution $v$ on any interval $(0,T)$.
\item For every $\kappa>1$,
 \begin{equation}\label{sandwich}
  \Ts(\kappa u_0)\ \le\ \lim_{n\to\infty}\Tmax(\min\{u_0,n\})\ \le\ \Ts(u_0).
 \end{equation}
 More generally, $\lim_n\Tmax(w_n)\le\Ts(u_0)$ for every nondecreasing
 sequence of bounded nonnegative functions $w_n$ converging to $u_0$
 a.e.
\item The equality $\lim_n\Tmax(\min\{u_0,n\})=\Ts(u_0)$ holds for every
 $0\le u_0\in L^1(\Omega)$ if and only if \eqref{CB} holds.
\item Let $\kappa_*(u_0):=\sup\{\kappa>0:\ \Ts(\kappa u_0)>0\}\in[0,\infty]$.
 Then $\kappa u_0$ satisfies \textnormal{(i)--(iii)} for every
 $\kappa>\kappa_*(u_0)$, and \textnormal{(i)} fails for $\kappa u_0$ for
 every $0<\kappa<\kappa_*(u_0)$.
\end{enumerate}
\end{maintheorem}

Part (a) says that $\Ts(u_0)$ is attained and sharply separates
convergence from divergence of the truncated solutions; for bounded data
it is the complete blow-up time of Martel \cite{Ma98}. Part (d) says that,
along every ray of data, the positive answer to Open Problem~5 occurs
exactly where local solvability fails, up to the single value
$\kappa=\kappa_*$.

The third result is a quantitative concentration estimate.

\begin{maintheorem}[Critical concentration estimate]\label{thm:conc}
There is a constant $C_d\ge1$, depending only on $d$, with the following
property. Let $0\le u_0\in L^1(\Omega)$, $x_0\in\Omega$,
$R_0:=\operatorname{dist}(x_0,\partial\Omega)$, $\rho>0$, and
$\mu:=\int_{B_\rho(x_0)}u_0>0$. If
\[
 \tau:=C_d\,\rho^2\exp\!\big(C_d\,\mu^{-2/d}\big)<R_0^2,
\]
then $\Ts(u_0)\le\tau$. In particular, $\Tmax(u_0)\le\tau$ if $u_0$ is
bounded, and $U\equiv+\infty$ on $(\tau,\infty)\times\Omega$.
Equivalently: if \eqref{P} has a nonnegative integral solution on
$(0,T)$ and $T':=\min\{T,R_0^2\}$, then
\begin{equation}\label{necessary}
 \int_{B_\rho(x_0)}u_0\le C_d^{d/2}
 \Big[\log\frac{T'}{C_d\rho^2}\Big]^{-d/2}
 \qquad\text{whenever }C_d\rho^2<T'.
\end{equation}
The exponent $-2/d$ cannot be improved.
\end{maintheorem}

Estimate \eqref{necessary} is the interior part of the necessary
condition of Hisa and Ishige \cite{HI18,HI24}; our proof is different.
For bounded data the bound does not involve $\|w_0\|_\infty$, only the
mass in a small ball.

\begin{maintheorem}[Criteria, rates and the Orlicz scale]\label{thm:data}
Let $x_0\in\Omega$, $0\le u_0\in L^1(\Omega)$, $m(r):=\int_{B_r(x_0)}u_0$
and $\mu_d:=(C_d/2)^{d/2}$.
\begin{enumerate}[label=\textnormal{(\alph*)}]
\item If $\limsup_{r\downarrow0}(\log\frac1r)^{d/2}m(r)>\mu_d$, then
 $\Ts(w_n)\to0$ for every sequence $0\le w_n\in L^1(\Omega)$ with
 $w_n\to u_0$ in $L^1(B_\delta(x_0))$ for some $\delta>0$. In particular
 $u_0$ satisfies \textnormal{(i)--(iii)}.
\item For $\gamma>1$ and $A>0$ let
 \[
  u_{\gamma,A}(x):=A\,|x-x_0|^{-d}
  \Big(\log\frac1{|x-x_0|}\Big)^{-\gamma}\ \text{ if }|x-x_0|<\delta,
  \qquad u_{\gamma,A}(x):=0\ \text{ otherwise},
 \]
 where $0<\delta\le\min\{e^{-1-\gamma/d},\operatorname{dist}(x_0,\partial\Omega)\}$.
 If $1<\gamma<1+d/2$, then $u_{\gamma,A}$ satisfies
 \textnormal{(i)--(iii)} for every $A>0$, and there are $C,n_0$ such that
 \begin{equation}\label{rate}
  \frac d2\,n^{-2/d}\le\Tmax(\min\{u_{\gamma,A},n\})
  \le n^{-2/d}\exp\!\Big(C(\log n)^{\frac{2\gamma}{d+2}}\Big)
  \qquad(n\ge n_0),
 \end{equation}
 where $\frac{2\gamma}{d+2}<1$. If $\gamma=1+d/2$, there is a critical
 coefficient $A_c\in[0,A_d]$, $A_d:=d\mu_d/(2|\mathbb S^{d-1}|)$, such
 that $u_{\gamma,A}$ satisfies \textnormal{(i)--(iii)} for $A>A_c$ and
 \textnormal{(i)} fails for $A<A_c$; moreover $A_c>0$ by
 \cite{HI18,HI24}.
\item For every $q\in[0,d/2)$ there is $0\le u_0\in L^1(\Omega)$ with
 $\int_\Omega u_0[\log(e+u_0)]^q<\infty$ satisfying
 \textnormal{(i)--(iii)}. On the other hand, if $u_0\ge0$ and
 $\int_\Omega u_0[\log(e+u_0)]^{d/2}<\infty$, then $\Ts(u_0)>0$, and
 \textnormal{(i)--(iii)} fail.
\end{enumerate}
\end{maintheorem}

By \eqref{rate}, the truncation lifespans of $u_{\gamma,A}$ collapse at
the fastest rate compatible with the height of the truncation, up to the
factor $\exp(C(\log n)^{2\gamma/(d+2)})=n^{o(1)}$. Part (c) places the
dividing line of Open Problem~5 exactly at $L\log^{d/2}L$, where Miyamoto
\cite{Mi21} proved local solvability.

Finally, Theorem~\ref{thm:conc} applied to vanishing concentrated
bumps gives ill-posedness at every datum.

\begin{maincorollary}[Ill-posedness at every datum]\label{cor:AR}
Let $\nu$ be a finite nonnegative Borel measure on $\Omega$. There are
$0\le w_n\in C_c^\infty(\Omega)$ such that $w_n\to\nu$ narrowly,
$w_n\to\nu$ in $L^1(\Omega)$ if $\nu\in L^1(\Omega)$, and
$\Ts(w_n)\le1/n$. In particular $\Tmax(w_n)\to0$, and the minimal
solutions with data $w_n$ are identically infinite on
$(1/n,\infty)\times\Omega$.
\end{maincorollary}

\subsection{Comments on the main results}

\begin{remark}[On Theorem~\ref{thm:main}]\label{rem:main-b}
(1) There is no intermediate regime in which $U$ is finite on part of
$\Omega$ only, or leaves $L^1_{\rm loc}(\Omega)$ and later returns. For
general $u_0\in L^1$, $\Ts(u_0)$ is the existence time of the proper
minimal solution in the sense of \cite[Definition~2.2]{AR04}.

(2) If $\kappa\mapsto\Ts(\kappa u_0)$ is right-continuous at $\kappa=1$,
then equality holds in \eqref{sandwich} without \eqref{CB}; the only
obstruction to equality is a jump of $\Ts(\kappa u_0)$ at $\kappa=1^+$.

(3) Monotonicity of the approximations matters in
Theorem~\ref{thm:main}(b): by Corollary~\ref{cor:AR}, every $u_0\in L^1$,
including $u_0=0$, is an $L^1$-limit of smooth data $w_n$ with
$\Tmax(w_n)\to0$.
\end{remark}

\begin{remark}[On Theorem~\ref{cor:OP5}]\label{rem:OP5}
(1) Nonexistence in the class of \cite[Theorem~11]{BC96} is a priori
weaker than (iii), since integral solutions need not be continuous with
values in $L^1$, nor bounded for positive times. The datum of
\cite{BC96} nevertheless satisfies (iii), because the divergence
\eqref{BC-claim} on which its construction rests also excludes integral
solutions.

(2) Property (iii) concerns nonnegative solutions. Whether a datum
satisfying (iii) admits a sign-changing solution of
$u_t-\Delta u=|u|^{p-1}u$ is a different question, related to Open
Problems~2 and 4 of \cite{BC96} (see \cite{To} for the failure of
positivity preservation). It plays no role in Open Problem~5.
\end{remark}

\begin{remark}[On Theorem~\ref{thm:conc}]\label{rem:conc}
(1) \emph{Scale invariance.} The map
$u\mapsto\lambda^du(\lambda^2t,\lambda x)$ preserves the equation and the
mass $\mu$, and multiplies $\tau$ and $R_0^2$ by $\lambda^{-2}$; this is why
$C_d$ does not depend on $\Omega$.

(2) \emph{Two regimes.} If $\mu$ is above a fixed threshold, then
$\tau\approx\rho^2$ and the mass blows up on its own parabolic time scale,
as in Kaplan's argument. If $\mu$ is small, the factor
$\exp(C_d\mu^{-2/d})$ accounts for the parabolic scales needed to reach
that threshold.

(3) \emph{Interior only.} The condition $\tau<R_0^2$ degenerates as
$x_0\to\partial\Omega$; data concentrating at the boundary obey a
different, weighted condition \cite{HI24}.
\end{remark}

\begin{remark}[On the hypothesis \eqref{CB}]\label{rem:CB}
Property \eqref{CB} is used only in Theorem~\ref{thm:main}(c), and hence in Theorem~\ref{cor:OP5}(c). If $u_0$
satisfies (i) but not (iii), choose $k$ with
$\Tmax(u_{0,k})<\Ts(u_0)$. Since $\Ts$ is nonincreasing with respect to
the datum (Lemma~\ref{lem:U} and Theorem~\ref{thm:main}(a)),
$\Ts(u_{0,k})\ge\Ts(u_0)>\Tmax(u_{0,k})$, so the bounded datum $u_{0,k}$
violates \eqref{CB}. Martel \cite{Ma98} proved complete blow-up in
general smooth bounded domains for solutions that are nondecreasing in
time, which the truncations $u_{0,k}$ need not produce. The results of
\cite{BCo87} are quoted for convex domains in \cite{Du19} and without
this restriction in \cite[p.~744]{AR04}, and complete blow-up of all
finite-time blow-up solutions for $1<p<p_S$ is stated in
\cite[Theorem~A]{Du19}. If
\eqref{CB} holds in every smooth bounded domain, then Theorem~\ref{thm:main}(c) and the equivalence in
Theorem~\ref{cor:OP5} hold without restriction.
\end{remark}

\subsection{Relation to the literature}\label{subsec:lit}
\label{subsec:AR}

\emph{Arrieta and Rodr\'{\i}guez-Bernal.} In \cite{AR04} the existence
time is shown not to be uniform on bounded subsets of critical spaces,
which answers Open Problem~1 of \cite{BC96}; in supercritical ranges
they obtain ill-posedness at every nonnegative datum
\cite[Theorem~1.1(i)]{AR04}. At $p=1+2/d$ they construct smooth data
converging to $\nu+C\delta_{x_0}$, with some $C>0$, whose existence times
tend to zero \cite[Theorem~3.9]{AR04}. Their solutions are the proper
minimal solutions, which coincide with our $U$ (Remark~\ref{rem:proper}),
and their main tool is a Kaplan-type lemma \cite[Lemma~3.3]{AR04}, of
which Lemma~\ref{lem:kaplan} is a scale-invariant form. That lemma
requires a mass above a fixed threshold in a ball, independent of its
radius; for a fixed $u_0\in L^1$ the mass of small balls tends to zero, so
it does not apply at small scales. Theorem~\ref{thm:conc} quantifies how a
small mass is amplified across parabolic scales until the threshold is
reached, and Corollary~\ref{cor:AR} removes the Dirac mass in
\cite[Theorem~3.9]{AR04}. On the other hand, \cite{AR04} treats general
nonlinearities $f(u)\ge au^p-bu$, unbounded domains, weighted spaces and
nonlinear boundary conditions.

\emph{Instantaneous blow-up.} The second formulation of Open Problem~5
has a linear prototype. Baras and Goldstein \cite{BG84} showed that for
$u_t-\Delta u=c|x|^{-2}u$ with $c>(d-2)^2/4$ the solutions with truncated
potential $\min\{c|x|^{-2},n\}$ tend to $+\infty$ at every point and every
positive time, which they called instantaneous blow-up. Cabr\'e and
Martel \cite{CM99} related existence versus instantaneous blow-up for
general potentials to a Hardy-type inequality. For nonlinear problems,
nonexistence and instantaneous blow-up for data above a singular steady
state were obtained by Peral and V\'azquez \cite{PV95} and V\'azquez
\cite{Va99} for $u_t-\Delta u=e^u$, and by Galaktionov and V\'azquez
\cite{GV97} for $u^p$ in supercritical ranges. There the obstruction is a
singular potential or a singular steady state; in Open Problem~5 it is
the singularity of the $L^1$ datum, and Theorem~\ref{thm:main}(a) is the
corresponding statement for $u^p$ at the doubly critical exponent.

\emph{Complete blow-up.} For convex nonlinearities and bounded data,
Martel \cite{Ma98}, building on \cite{BCMR96}, showed that the limit of
the truncated problems is a weak solution as long as blow-up is not
complete, and that no continuation exists after the complete blow-up
time. Theorem~\ref{thm:main}(a) is an analogue for $L^1$ data and
integral solutions. Complete blow-up after $\Tmax$ for bounded data and
$1<p<p_S$ was proved by Baras and Cohen \cite{BCo87}; see also
\cite{GV97} and \cite[Section~27]{QS19}.

\emph{Necessary conditions and optimal singularities.} Necessary
conditions for nonnegative solutions were given by Baras and Pierre
\cite{BP85} in terms of capacities, and by Hisa and Ishige
\cite{HI18,HI24} in the form \eqref{necessary}, for the Cauchy problem
and for smooth domains; they also identify the optimal singularity
$|x|^{-d}(\log\frac1{|x|})^{-1-d/2}$. Miyamoto \cite{Mi21} proved
solvability in $L\log^{d/2}L$ and nonexistence for some data in
$L\log^qL$, $q<d/2$; see also \cite{FI18}. Near $\partial\Omega$ the
admissible singularities differ from those in the interior \cite{HI24}.
Laister, Robinson,
Sier\.{z}\c{e}ga and Vidal-L\'opez \cite{LRSV16} characterised local
existence in Lebesgue spaces, using data close to those of
Remark~\ref{rem:spikes}; see also \cite{LS20}.

To the best of our knowledge, the inequalities \eqref{sandwich}, and
hence the fact that nonexistence of integral solutions yields a positive
answer to Open Problem~5, have not been recorded before.

\subsection{Strategy of the proofs}\label{sec:strategy}

The proofs rest on four ideas.

\emph{1. No return.} If the minimal solution $U$ is finite at one point
$(t_0,x_0)$, then, by the Chapman--Kolmogorov identity,
$\int_\Omega G(t_0-t,x_0,y)U(t,y)\dd y\le U(t_0,x_0)$ for $t<t_0$, so $U$
is locally integrable, hence an integral solution, up to $t_0$. Thus
$U$ is finite before $\Ts(u_0)$ and infinite everywhere afterwards
(Theorem~\ref{thm:main}(a)), which gives
(ii)$\Leftrightarrow$(iii).

\emph{2. A data-independent bound.} A bounded solution that exists up
to time $t$ satisfies $S(t)w_0\le(2t/d)^{-d/2}$ (Weissler). Applied to the
solutions with data $\min\{u_0,n\}$, it bounds their weighted
$L^1$ norms uniformly in $n$ on a common existence interval, so that a
common interval produces an integral solution in the limit. This is the
right inequality in \eqref{sandwich} and the implication
(iii)$\Rightarrow$(i).

\emph{3. A concave transformation.} The function
$\Phi(s)=(s^{1-p}+c)^{-1/(p-1)}$ satisfies $\Phi'(s)s^p=\Phi(s)^p$ and is
bounded. Applied to $U$, it turns an integral solution, which may be
unbounded, into a bounded supersolution with datum $\Phi(u_0)\ge\lambda
\min\{u_0,n\}$ (Lemma~\ref{lem:concave}). This gives the left inequality
in \eqref{sandwich}, and is used again to pass from bounded to general
data in Theorem~\ref{thm:conc}.

\emph{4. Iteration over parabolic scales.} At $p=1+2/d$ the Duhamel term
gains a factor $\int_{\rho^2}^{h}s^{d/2-dp/2}\dd s=\log(h/\rho^2)$: a mass
$\mu$ in $B_\rho$ produces, after a time $h\approx\rho^2e^{2/(\beta\mu^{p-1})}$,
a mass $2\mu$ in $B_{b\sqrt h}$ (Lemma~\ref{lem:onestep}). Iterating about
$\log_2(1/\mu)$ times until the mass exceeds the threshold of a local
Kaplan criterion (Lemma~\ref{lem:kaplan}) costs a total time of order
$\rho^2\exp(C\mu^{-2/d})$. This is Theorem~\ref{thm:conc}; the criteria and
rates of Theorem~\ref{thm:data} follow by choosing $\rho$ optimally.

The logical structure is
\[
 \text{(ii)}\ \overset{1}{\Longleftrightarrow}\ \text{(iii)}\
 \overset{2}{\Longrightarrow}\ \text{(i)}\ \overset{3}{\Longrightarrow}\
 \Ts(\kappa u_0)=0\ \text{ for all }\kappa>1,
 \qquad
 \text{(i)}\ \overset{\eqref{CB}}{\Longrightarrow}\ \text{(iii)} .
\]

\section{Preliminaries}\label{sec:prelim}

\subsection{Notation and constants}\label{subsec:constants}

Throughout, $G=G_\Omega$ is the Dirichlet heat kernel and $S(t)$ the
Dirichlet heat semigroup; $\lambda_1$ and $\varphi_1>0$ are the first
Dirichlet eigenvalue and eigenfunction, with $\int_\Omega\varphi_1=1$;
$\omega_d:=|B_1|$ and $\sigma_d:=|\mathbb S^{d-1}|$; $R_0$ denotes the
distance of a point $x_0$ to $\partial\Omega$. For a datum $u_0$ we write
$u_{0,n}=\min\{u_0,n\}$, $U_n$ for the solutions of \eqref{Pn} and
$U=\lim_nU_n$ for the minimal solution (Lemma~\ref{lem:U}). The following
constants depend only on $d$ and are fixed throughout:
\begin{center}
\small
\begin{tabular}{@{}lll@{}}
\toprule
constant & definition & used in\\
\midrule
$b,\ \Lambda,\ M_d$ & $b=2^{-7}$, $\Lambda=\lambda_1(B_1)/b^2$,
 $M_d=\max\{(2\Lambda)^{d/2},(2d)^{d/2}\}$ & Lemma~\ref{lem:kaplan}\\
$a,\ c_0,\ c_1,\ c_2$ & $a=2^{-4}$, $c_0=(4\pi)^{-d/2}e^{-d^2\pi^2/64}$,
 $c_1=c_0e^{-1/4}$, $c_2=c_1^{1+p}\omega_da^d$ & Lemma~\ref{lem:onestep}\\
$\beta,\ \mu_1$ & $\beta=\omega_db^dc_2$, $\mu_1=M_d/c_1$ &
 proof of Theorem~\ref{thm:conc}\\
$C_d$ & see \eqref{CN} & Theorems~\ref{thm:conc}, \ref{thm:data}\\
$\mu_d,\ A_d$ & $\mu_d=(C_d/2)^{d/2}$, $A_d=d\mu_d/(2\sigma_d)$ &
 Theorem~\ref{thm:data}\\
\bottomrule
\end{tabular}
\end{center}
where $\lambda_1(B_1)$ is the first Dirichlet eigenvalue of the unit
ball, and
\begin{equation}\label{CN}
 C_1:=\frac{2}{\beta\,(1-2^{-2/d})},\qquad
 C_2:=\frac{d\mu_1^{2/d}}{2\log2},\qquad
 C_d:=2^{10}\max\{48,\,C_1+C_2\}.
\end{equation}

\subsection{The Dirichlet heat semigroup}

Let $G=G_\Omega(t,x,y)$ be the Dirichlet heat kernel of $\Omega$. It is
smooth and positive on $(0,\infty)\times\Omega\times\Omega$, symmetric
in $(x,y)$, and satisfies the Chapman--Kolmogorov identity
\begin{equation}\label{CK}
 \int_\Omega G(t,x,z)\,G(s,z,y)\dd z=G(t+s,x,y)
 \qquad(t,s>0,\ x,y\in\Omega),
\end{equation}
which is the semigroup property $S(t)S(s)=S(t+s)$ written for the
kernel; see e.g.\ \cite[Chapters~2 and 3]{Da89}. Moreover
$\int_\Omega G(t,x,y)\dd y\le1$ and
$G(t,x,y)\le(4\pi t)^{-d/2}e^{-|x-y|^2/(4t)}$. For measurable $f\ge0$
we set
\[
 S(t)f(x):=\int_\Omega G(t,x,y)f(y)\dd y\in[0,\infty];
\]
for $f\in L^1(\Omega)$ this is the Dirichlet heat semigroup, and
$\|S(t)f\|_\infty\le(4\pi t)^{-d/2}\|f\|_1$. Since
$\int_\Omega G(t,x,y)\dd y\le1$, H\"older's inequality gives
\begin{equation}\label{jensen}
 \big([S(t)f](x)\big)^p\le[S(t)f^p](x)\qquad\text{for all measurable }f\ge0 .
\end{equation}

Let $\lambda_1>0$ be the first Dirichlet eigenvalue and $\varphi_1>0$
the corresponding eigenfunction, normalised by
$\int_\Omega\varphi_1=1$. Since $S(t)\varphi_1=e^{-\lambda_1t}\varphi_1$
and $G$ is symmetric, Tonelli's theorem gives
\begin{equation}\label{phi1}
 \int_\Omega\big(S(t)f\big)\varphi_1=e^{-\lambda_1t}\int_\Omega f\varphi_1
 \qquad\text{for all measurable }f\ge0,\ t\ge0 .
\end{equation}

We use the following explicit interior lower bound, due to van den Berg
\cite{vdB89,vdB90}; see \cite[(6)]{LRSV16}: if the segment joining
$x,y\in\Omega$ has distance at least $\ell>0$ from $\partial\Omega$, then
\begin{equation}\label{vdB}
 G(t,x,y)\ge e^{-d^2\pi^2t/(4\ell^2)}(4\pi t)^{-d/2}e^{-|x-y|^2/(4t)}
 \qquad\text{for all }t>0.
\end{equation}
\subsection{Bounded data}

We use the following standard facts; see e.g.\ \cite[Chapter~II]{QS19}
and \cite[Section~1]{BC96}. Item (d) is the comparison principle, and
(f) follows by comparison with the spatially constant solution
$\big(n^{1-p}-(p-1)t\big)^{-1/(p-1)}$, $n:=\|w_0\|_\infty$, which exists
up to time $n^{1-p}/(p-1)=\frac d2n^{-2/d}$. Item (g) is Weissler's
necessary condition \cite{We81}. Its proof uses only the semigroup
property and \eqref{jensen}, so it applies verbatim to the Dirichlet
semigroup: the function $Y(s):=[S(t-s)w(s)](x)$ satisfies
$Y(s)\ge Y(0)+\int_0^sY(\sigma)^p\dd\sigma$ on $[0,t]$, and integrating
this ODE inequality gives $Y(0)^{1-p}\ge(p-1)t$.

\begin{lemma}\label{lem:bounded}
Let $0\le w_0\in L^\infty(\Omega)$. There is a unique maximal
nonnegative solution $w$ of \eqref{P} on $[0,\Tmax(w_0))$, and:
\begin{enumerate}[label=\textnormal{(\alph*)}]
\item $w\in L^\infty(\Omega\times(0,T'))$ for every $T'<\Tmax(w_0)$;
 $w\in C^{1,2}((0,\Tmax)\times\Omega)\cap C((0,\Tmax)\times\overline\Omega)$
 and $w=0$ on $\partial\Omega$;
\item for $0\le s<t<\Tmax(w_0)$ and $x\in\Omega$, with $w(0):=w_0$,
 \begin{equation}\label{mild}
  w(t,x)=\big[S(t-s)w(s)\big](x)+\int_s^t\big[S(t-\sigma)w(\sigma)^p\big](x)\dd\sigma;
 \end{equation}
\item if $\Tmax(w_0)<\infty$, then $\|w(t)\|_\infty\to\infty$ as
 $t\uparrow\Tmax(w_0)$;
\item if $0\le a_0\le b_0$ are bounded, the corresponding solutions
 satisfy $a\le b$ on their common interval of existence, and
 $\Tmax(b_0)\le\Tmax(a_0)$;
\item for $t_1\in(0,\Tmax(w_0))$, $w(t_1+\cdot)$ is the solution with
 datum $w(t_1)$, and $\Tmax(w(t_1))=\Tmax(w_0)-t_1$;
\item if $w_0\not\equiv0$, then
 \[
  \Tmax(w_0)\ge\frac d2\,\|w_0\|_\infty^{-2/d};
 \]
\item for $0<t<\Tmax(w_0)$,
 \begin{equation}\label{weissler}
  S(t)w_0\le\big((p-1)t\big)^{-\frac1{p-1}}=\Big(\frac{2t}d\Big)^{-d/2}
  \qquad\text{in }\Omega .
 \end{equation}
\end{enumerate}
\end{lemma}

\begin{remark}\label{rem:classes}
(1) By \eqref{mild} with $s=0$, the solution in Lemma~\ref{lem:bounded}
is an integral solution on $(0,\Tmax(w_0))$; hence
$\Tmax(w_0)\le\Ts(w_0)$.

(2) Let $u\ge0$ belong to
$C([0,T];L^1(\Omega))\cap L^\infty_{\rm loc}((0,T);L^\infty(\Omega))$ and
be a classical solution of \eqref{P} on $(0,T)$, with $u(0)=u_0$. For
$0<s<t<T$, \eqref{mild} holds by uniqueness of bounded solutions. As
$s\downarrow0$, $S(t-s)u(s)\to S(t)u_0$ in $L^1$, while
$\int_s^tS(t-\sigma)u(\sigma)^p\dd\sigma$ increases to
$\int_0^tS(t-\sigma)u(\sigma)^p\dd\sigma$. Hence $u$ is an integral
solution on $(0,T)$. Consequently, property (iii) is stronger than
nonexistence in the class of \cite[Theorem~11]{BC96}.
\end{remark}

\subsection{The truncated problems and the minimal solution}

\begin{lemma}\label{lem:Un}
Let $0\le u_0\in L^1(\Omega)$ and $n\in\N$. There is a unique
nonnegative function $U_n$ in $C([0,\infty);L^1(\Omega))$ such that
\begin{equation}\label{Un-mild}
 U_n(t)=S(t)u_0+\int_0^tS(t-s)g_n(U_n(s))\dd s\qquad(t\ge0).
\end{equation}
Moreover:
\begin{enumerate}[label=\textnormal{(\alph*)}]
\item the Picard iterates $P_0:=0$,
 $P_{k+1}(t):=S(t)u_0+\int_0^tS(t-s)g_n(P_k(s))\dd s$ are nonnegative,
 nondecreasing in $k$, and converge to $U_n$ in $C([0,T];L^1(\Omega))$
 for every $T>0$;
\item $U_n$ is nondecreasing with respect to $u_0$ and with respect to
 $n$;
\item the right-hand side of \eqref{Un-mild}, evaluated pointwise, is a
 representative of $U_n$ that is continuous on $(0,\infty)\times\Omega$.
 We always use this representative;
\item $U_n$ is bounded on $(t_1,t_2)\times\Omega$ for $0<t_1<t_2$,
 $U_n\in C^{2,1}((0,\infty)\times\Omega)\cap C((0,\infty)\times\overline\Omega)$,
 $U_n=0$ on $(0,\infty)\times\partial\Omega$, and $U_n>0$ in
 $(0,\infty)\times\Omega$ unless $u_0=0$ a.e.
\end{enumerate}
\end{lemma}

Lemma~\ref{lem:Un} is standard. Since $g_n$ is bounded and globally
Lipschitz, existence, uniqueness and the convergence in (a) follow from
the Banach fixed point theorem in $C([0,T];L^1(\Omega))$ with the norm
$\sup_{t\le T}e^{-2Lt}\|u(t)\|_1$, where $L=pn^{(p-1)/p}$ is the
Lipschitz constant of $g_n$. The monotonicity in (a) and (b) follows by
induction on the Picard iterates, since $S(t)$ is order preserving and
$g_n\le g_{n+1}$ are nondecreasing. For (c), the Duhamel term is the heat
potential of a bounded function. For (d), $0\le U_n(t)\le(4\pi
t)^{-d/2}\|u_0\|_1+nt$, and for $t>t_1$ the function $U_n(t_1+\cdot)$
solves \eqref{Pn} with a bounded datum and a bounded Lipschitz
nonlinearity, so it is classical by standard parabolic regularity;
finally $U_n(t)\ge S(t)u_0>0$ unless $u_0=0$ a.e.

\begin{lemma}[Comparison]\label{lem:dominate}
Let $T\in(0,\infty]$ and let $v:\Omega\times(0,T)\to[0,\infty]$ be
measurable with
\begin{equation}\label{super}
 v(t)\ge S(t)u_0+\int_0^tS(t-s)v(s)^p\dd s\qquad\text{a.e. in }\Omega\times(0,T).
\end{equation}
Then $U_n\le v$ a.e. in $\Omega\times(0,T)$ for every $n$. In particular
this holds for every nonnegative integral solution $v$ on $(0,T)$.
\end{lemma}

\begin{proof}
Let $\bar v$ be the right-hand side of \eqref{super}, evaluated pointwise
in $[0,\infty]$; then $\bar v\le v$ a.e. We show $P_k\le v$ a.e. by
induction. This is clear for $k=0$. If $P_k\le v$ a.e., then
$g_n(P_k)\le P_k^p\le v^p$ a.e., so $P_{k+1}\le\bar v$ everywhere and
$P_{k+1}\le v$ a.e. Letting $k\to\infty$ and using
Lemma~\ref{lem:Un}(a) gives the claim.
\end{proof}

\begin{lemma}\label{lem:U}
Let $0\le u_0\in L^1(\Omega)$. By Lemma~\ref{lem:Un}(b),(c) the limit
$U:=\lim_{n\to\infty}U_n\in[0,\infty]$ exists at every point of
$(0,\infty)\times\Omega$. It is measurable, nondecreasing with respect to
$u_0$, and
\begin{equation}\label{U-eq}
 U(t,x)=[S(t)u_0](x)+\int_0^t[S(t-s)U(s)^p](x)\dd s\qquad
 \text{for all }(t,x)\in(0,\infty)\times\Omega .
\end{equation}
In particular, if $U<\infty$ a.e. in $\Omega\times(0,T)$, then $U$ is a
nonnegative integral solution of \eqref{P} on $(0,T)$. Moreover, $U\le v$
a.e. for every $v$ as in Lemma~\ref{lem:dominate}.
\end{lemma}

\begin{proof}
The functions $g_n(U_n)$ are nondecreasing in $n$ and converge pointwise
to $U^p$ (with $\infty^p:=\infty$): if $U(s,y)<\infty$, then
$\min\{U_n(s,y)^p,n\}\to U(s,y)^p$, and if $U(s,y)=\infty$, then
$\min\{U_n(s,y)^p,n\}\to\infty$. Monotone convergence in \eqref{Un-mild},
written through the kernel $G$, gives \eqref{U-eq}. The last assertion
follows from Lemma~\ref{lem:dominate}.
\end{proof}

\begin{remark}[Proper minimal solution]\label{rem:proper}
The limit $U$ does not depend on the truncation. Let
$f_n:[0,\infty)\to[0,\infty)$ be nondecreasing and globally Lipschitz,
with $f_n\le f_{n+1}\le s^p$ and $f_n(s)\to s^p$ for every $s\ge0$; for
instance, $f_n(s)=s^p$ for $s\le n$ and
$f_n(s)=n^p+pn^{p-1}(s-n)$ for $s>n$, as in \cite{BCo87,AR04}. Let $V_n$
be the corresponding solutions and $V:=\lim_nV_n$. The Picard iterates
for $f_n$ stay below $U$ by \eqref{U-eq}, since $f_n(P_k)\le P_k^p\le
U^p$; hence $V\le U$. Conversely, for $n\ge m$ we have
$f_m(V_n)\le f_n(V_n)\le V^p$; letting first $n\to\infty$ and then
$m\to\infty$ gives $f_n(V_n)\to V^p$ in $[0,\infty]$. By monotone
convergence $V$ satisfies \eqref{U-eq}, hence $U\le V$ by
Lemma~\ref{lem:dominate}.
Thus $V=U$, and $U$ is the proper minimal solution of
\cite{BCo87,GV97,AR04}.
\end{remark}

\begin{remark}[Bounded data]\label{rem:bdd-U}
Let $0\le w_0\in L^\infty(\Omega)$ with solution $w$. Then $U=w$ on
$(0,\Tmax(w_0))\times\Omega$. Indeed, $U\le w$ a.e. by
Lemma~\ref{lem:U} and Remark~\ref{rem:classes}(1). Conversely, for
$T'<\Tmax(w_0)$ and $n\ge\sup_{[0,T']}\|w\|_\infty^p$ we have
$g_n(w)=w^p$ on $[0,T']$, so $U_n=w$ there by the uniqueness in
Lemma~\ref{lem:Un}; hence $U\ge w$, and $U=w$ by continuity.
\end{remark}

\section{Proof of Theorems~\ref{cor:OP5} and \ref{thm:main}}\label{sec:main}

\subsection{Proof of Theorem~\ref{thm:main}(a)}

\begin{proof}[Proof of Theorem~\ref{thm:main}\textnormal{(a)}]
\emph{Step 1: no return.} Suppose that $U(t_0,x_0)<\infty$ for some
$t_0>0$ and $x_0\in\Omega$, and let $0<t<t_0$. By \eqref{CK} and Tonelli's
theorem, pointwise in $\Omega$,
\[
 S(t_0)u_0=S(t_0-t)\big[S(t)u_0\big],\qquad
 S(t_0-s)U(s)^p=S(t_0-t)\big[S(t-s)U(s)^p\big]\quad(0<s<t).
\]
Splitting the time integral in \eqref{U-eq} at $s=t$ and using
\eqref{U-eq} at time $t$, we obtain
\begin{equation}\label{noreturn}
\begin{aligned}
 U(t_0,x_0)
 &=\Big[S(t_0-t)\Big(S(t)u_0+\int_0^tS(t-s)U(s)^p\dd s\Big)\Big](x_0)\\
 &\qquad+\int_t^{t_0}\big[S(t_0-s)U(s)^p\big](x_0)\dd s\\
 &\ge\big[S(t_0-t)U(t)\big](x_0)
 =\int_\Omega G(t_0-t,x_0,y)\,U(t,y)\dd y .
\end{aligned}
\end{equation}
Since $G(t_0-t,x_0,\cdot)$ is continuous and positive in $\Omega$, it is
bounded below by a positive constant on every compact subset of
$\Omega$, so $U(t)\in L^1_{\rm loc}(\Omega)$ for every $t\in(0,t_0)$. In
particular $U<\infty$ a.e. in $\Omega\times(0,t_0)$, and by
Lemma~\ref{lem:U}, $U$ is an integral solution on $(0,t_0)$.

\emph{Step 2: conclusion.} Let $0<T<\Ts(u_0)$. By definition there is an
integral solution $v$ on $(0,T')$ for some $T'>T$, and $U\le v<\infty$
a.e. in $\Omega\times(0,T')$ by Lemma~\ref{lem:U}. Choosing $t_0\in(T,T')$
and $x_0$ with $U(t_0,x_0)<\infty$, Step~1 shows that
$U(t)\in L^1_{\rm loc}(\Omega)$ for $t\le T$ and that $U$ is an integral
solution on $(0,T)$. As $T<\Ts(u_0)$ is arbitrary, $U$ is an integral
solution on $(0,\Ts(u_0))$. Minimality is Lemma~\ref{lem:U}. Finally, if
$U(t_0,x_0)<\infty$ for some $t_0>\Ts(u_0)$, then $U$ is an integral
solution on $(0,t_0)$ by Step~1, which contradicts the definition of
$\Ts(u_0)$.
\end{proof}

\subsection{The concave transformation}

The following lemma is a parabolic form of the concave transformation of
\cite{BCMR96}. It is used in the proofs of both Theorem~\ref{thm:main}(b)
and Theorem~\ref{thm:conc}. The idea is as follows. The function $\Phi$
below is the backward flow of the ODE $y'=y^p$ over the time
$c/(p-1)$: if $y(0)=s$, then $y(-c/(p-1))=\Phi(s)$. Such a map commutes
with the reaction term, $\Phi'(s)s^p=\Phi(s)^p$; being concave, it turns
solutions into supersolutions; and being bounded by $c^{-1/(p-1)}$, it
turns a possibly unbounded integral solution into a bounded supersolution.
The price is that its initial value $\Phi(u_0)$ is only comparable with
$\lambda u_0$ for some $\lambda<1$.

\begin{lemma}[Concave transformation]\label{lem:concave}
Let $0\le u_0\in L^1(\Omega)$. Then $\Tmax(w)\ge\Ts(u_0)$ for every
$0\le w\in L^\infty(\Omega)$ such that $w\le\lambda u_0$ for some
$\lambda\in(0,1)$.
\end{lemma}

\begin{proof}
We may assume $\Ts(u_0)>0$ and $u_0\not\equiv0$. By
Theorem~\ref{thm:main}(a), $U$ is an integral solution on
$(0,T)$, $T:=\Ts(u_0)$. Fix $n\in\N$ with $n\ge\|w\|_\infty/\lambda$, so
that $w\le\lambda\min\{u_0,n\}$. Put $c:=(\lambda^{1-p}-1)n^{1-p}>0$,
$A:=c^{-1/(p-1)}$ and
\[
 \Phi(s):=\big(s^{1-p}+c\big)^{-\frac1{p-1}}
 =s\big(1+cs^{p-1}\big)^{-\frac1{p-1}}\quad(0\le s<\infty),\qquad
 \Phi(\infty):=A .
\]
Then $\Phi$ is continuous and nondecreasing on $[0,\infty]$, smooth on
$(0,\infty)$, and $\Phi(s)\le\min\{s,A\}$. A direct computation gives
\[
 \Phi'(s)=\big(1+cs^{p-1}\big)^{-\frac p{p-1}}\in(0,1],\qquad
 \Phi''(s)=-pc\,s^{p-2}\big(1+cs^{p-1}\big)^{-\frac{2p-1}{p-1}}<0,
\]
so $\Phi$ is concave and $1$-Lipschitz. Moreover
\begin{equation}\label{concave-id}
 \Phi'(s)\,s^p=s^p\big(1+cs^{p-1}\big)^{-\frac p{p-1}}
 =\Big[s\big(1+cs^{p-1}\big)^{-\frac1{p-1}}\Big]^p=\Phi(s)^p
 \qquad(s\ge0).
\end{equation}
Finally, $\Phi(s)\ge\lambda\min\{s,n\}$. Indeed, by the choice of $c$,
\[
 1+cn^{p-1}=1+\big(\lambda^{1-p}-1\big)=\lambda^{1-p},
\]
and since $s\mapsto(1+cs^{p-1})^{-1/(p-1)}$ is nonincreasing,
\[
 \Phi(s)\ge s\big(1+cn^{p-1}\big)^{-\frac1{p-1}}=\lambda s\quad(0\le s\le n),
 \qquad
 \Phi(s)\ge\Phi(n)=\lambda n\quad(s\ge n).
\]
Consequently $w\le\lambda\min\{u_0,n\}\le\Phi(u_0)$.

\emph{Step 1.} Let $U_k$ be the solutions of \eqref{Pn} (with $k$ in
place of $n$) and put $W_k:=\Phi(U_k)$ and
$F_k:=\Phi'(U_k)g_k(U_k)$. By \eqref{concave-id},
\[
 F_k=\Phi'(U_k)\min\{U_k^p,k\}
 =\min\big\{\Phi'(U_k)U_k^p,\,k\Phi'(U_k)\big\}
 =\min\big\{W_k^p,\,k\Phi'(U_k)\big\}.
\] We claim that
\begin{equation}\label{Wk-super}
 W_k(t,x)\ge[S(t)\Phi(u_0)](x)+\int_0^t[S(t-s)F_k(s)](x)\dd s
 \qquad\text{for all }(t,x)\in(0,\infty)\times\Omega .
\end{equation}
Fix $0<t_1<t$. By Lemma~\ref{lem:Un}(d), $U_k$ is a positive classical
solution in $(t_1,\infty)\times\Omega$, continuous on
$[t_1,\infty)\times\overline\Omega$ and zero on $\partial\Omega$. Since
$\Phi$ is smooth and concave on $(0,\infty)$,
\[
 \partial_tW_k-\Delta W_k=\Phi'(U_k)\big(\partial_tU_k-\Delta U_k\big)
 -\Phi''(U_k)|\nabla U_k|^2\ge\Phi'(U_k)g_k(U_k)=F_k\qquad\text{in }(t_1,\infty)\times\Omega,
\]
and $W_k=0$ on $\partial\Omega$. For $t\ge t_1$ let
\[
 Z(t):=W_k(t)-S(t-t_1)W_k(t_1)-\int_{t_1}^tS(t-s)F_k(s)\dd s .
\]
Since $F_k$ is bounded and continuous, $Z$ is continuous on
$[t_1,\infty)\times\overline\Omega$, $Z(t_1)=0$, $Z=0$ on
$\partial\Omega$, and $Z_t-\Delta Z\ge0$ in the sense of distributions in
$(t_1,\infty)\times\Omega$. By the weak maximum principle $Z\ge0$, that
is,
\[
 W_k(t)\ge S(t-t_1)W_k(t_1)+\int_{t_1}^tS(t-s)F_k(s)\dd s
 \qquad(t\ge t_1).
\]
Let $t_1\downarrow0$. Since $\Phi$ is $1$-Lipschitz and
$U_k\in C([0,\infty);L^1(\Omega))$, $W_k(t_1)\to\Phi(u_0)$ in $L^1$, so
$S(t-t_1)W_k(t_1)\to S(t)\Phi(u_0)$ pointwise in $\Omega$; the integral
converges by monotone convergence. This proves \eqref{Wk-super}.

\emph{Step 2.} Let $W:=\Phi(U)\in[0,A]$. Since $\Phi$ is continuous on
$[0,\infty]$, $W_k\to W$ everywhere. Since $U<\infty$ a.e. in
$\Omega\times(0,T)$, at almost every point of $\Omega\times(0,T)$ we have
$\Phi'(U_k)\to\Phi'(U)>0$, hence $k\Phi'(U_k)\to\infty$ and
$F_k\to W^p$. Since $0\le F_k\le A^p$, dominated convergence in
\eqref{Wk-super} yields
\begin{equation}\label{W-super}
 W(t,x)\ge[S(t)\Phi(u_0)](x)+\int_0^t[S(t-s)W(s)^p](x)\dd s
 \qquad\text{for all }(t,x)\in(0,T)\times\Omega .
\end{equation}

\emph{Step 3.} Let $z$ be the maximal classical solution with datum $w$
and $T_1:=\min\{T,\Tmax(w)\}$. Define $z_0:=0$ and
$z_{j+1}(t):=S(t)w+\int_0^tS(t-s)z_j(s)^p\dd s$. Since $w\le\Phi(u_0)$,
\eqref{W-super} gives $z_j\le W\le A$ on $(0,T)\times\Omega$ by
induction, and $z_j\le z$ on $(0,\Tmax(w))\times\Omega$ since $z$ solves
\eqref{mild}. The limit $z_\infty$ of the nondecreasing sequence $z_j$
is a solution of \eqref{mild} on $(0,T_1)$ with $0\le z_\infty\le A$. On
$[0,T']$ with $T'<T_1$, both $z$ and $z_\infty$ are bounded solutions of
\eqref{mild} with the same datum, so they coincide by Gronwall's
inequality. Hence $z\le A$ on $(0,T_1)\times\Omega$. If $\Tmax(w)<T$,
this contradicts Lemma~\ref{lem:bounded}(c). Therefore $\Tmax(w)\ge T$.
\end{proof}

\subsection{Proof of Theorem~\ref{thm:main}(b)--(d)}

\begin{proof}[Proof of Theorem~\ref{thm:main}\textnormal{(b)}]
\emph{Step 1: monotone approximations.} Let $0\le w_n\in L^\infty(\Omega)$
be nondecreasing in $n$ with $w_n\to u_0$ a.e., and let $w^{(n)}$ be the
corresponding solutions. By Lemma~\ref{lem:bounded}(d), $\Tmax(w_n)$ is
nonincreasing; let $T_\infty:=\lim_n\Tmax(w_n)$ and assume $T_\infty>0$.
Every $w^{(n)}$ is defined on $[0,T_\infty)$, and $w^{(n)}\le w^{(n+1)}$
there; let $u:=\lim_nw^{(n)}\in[0,\infty]$. Fix $0<t<T_\infty$ and put
$\tau:=(T_\infty-t)/2$. By Lemma~\ref{lem:bounded}(e), $w^{(n)}(t)$ is a
bounded datum with $\Tmax(w^{(n)}(t))=\Tmax(w_n)-t>\tau$, so
\eqref{weissler} gives $\|S(\tau)w^{(n)}(t)\|_\infty\le(2\tau/d)^{-d/2}$.
By \eqref{phi1} and $\int_\Omega\varphi_1=1$,
\[
 \int_\Omega w^{(n)}(t)\varphi_1
 =e^{\lambda_1\tau}\int_\Omega\big(S(\tau)w^{(n)}(t)\big)\varphi_1
 \le e^{\lambda_1\tau}\Big(\frac{2\tau}d\Big)^{-d/2},
\]
uniformly in $n$. By monotone convergence the same bound holds for $u(t)$,
and since $\varphi_1>0$, $u(t)<\infty$ a.e. in $\Omega$. Writing
\eqref{mild} with $s=0$ through the kernel $G$ and letting $n\to\infty$ by
monotone convergence in both terms, we see that $u$ satisfies
\eqref{integral} at every point of $\Omega\times(0,T_\infty)$. Hence $u$
is an integral solution on $(0,T_\infty)$, and $T_\infty\le\Ts(u_0)$. This
is the second assertion of (b); with $w_n=\min\{u_0,n\}$ it gives the
right inequality in \eqref{sandwich}.

\emph{Step 2: the left inequality.} Let $\kappa>1$ and apply
Lemma~\ref{lem:concave} to the datum $\kappa u_0$, with
$\lambda:=1/\kappa$ and $w:=\min\{u_0,n\}\le\lambda(\kappa u_0)$. This
gives $\Tmax(\min\{u_0,n\})\ge\Ts(\kappa u_0)$ for every $n$.
\end{proof}

\begin{proof}[Proof of Theorem~\ref{thm:main}\textnormal{(c)}]
Assume \eqref{CB}, and let $k\in\N$. Since $u_{0,k}\le u_0$,
Lemma~\ref{lem:U} gives $U\ge U^{(k)}$, where $U^{(k)}$ is the minimal
solution with datum $u_{0,k}$. By Theorem~\ref{thm:main}(a) and
\eqref{CB}, $U^{(k)}\equiv\infty$ on
$(\Tmax(u_{0,k}),\infty)\times\Omega$, hence so is $U$, and
Theorem~\ref{thm:main}(a) gives $\Ts(u_0)\le\Tmax(u_{0,k})$. Letting
$k\to\infty$ and using \eqref{sandwich} gives the equality. Conversely,
if the equality holds for every $u_0$, then for bounded $w_0$ we have
$\min\{w_0,n\}=w_0$ for large $n$, hence $\Tmax(w_0)=\Ts(w_0)$, which is
\eqref{CB}. The last assertion follows from Remark~\ref{rem:BCo}
below.
\end{proof}

\begin{proof}[Proof of Theorem~\ref{thm:main}\textnormal{(d)}]
By Lemma~\ref{lem:U} and part (a), $\Ts$ is nonincreasing with respect to
the datum, so $\kappa\mapsto\Ts(\kappa u_0)$ is nonincreasing. If
$\kappa>\kappa_*(u_0)$, then $\Ts(\kappa u_0)=0$, so $\kappa u_0$
satisfies (iii); it satisfies (ii) by part (a), and (i) by the right
inequality in \eqref{sandwich}. If
$0<\kappa<\kappa_*(u_0)$, choose $\kappa'>\kappa$ with
$\Ts(\kappa'u_0)>0$. The left inequality in \eqref{sandwich}, applied to
the datum $\kappa u_0$ with the ratio $\kappa'/\kappa>1$, gives
$\lim_n\Tmax(\min\{\kappa u_0,n\})\ge\Ts(\kappa'u_0)>0$, so (i) fails for
$\kappa u_0$.
\end{proof}

\begin{remark}[Baras--Cohen]\label{rem:BCo}
By Theorem~\ref{thm:main}(a) and $\Tmax(w_0)\le\Ts(w_0)$, property
\eqref{CB} is equivalent to complete blow-up in the sense of Baras and
Cohen: for every $0\le w_0\in L^\infty(\Omega)$ with $\Tmax(w_0)<\infty$,
the solutions $W_n$ of \eqref{Pn} with datum $w_0$ satisfy
$W_n(t,x)\to\infty$ for all $t>\Tmax(w_0)$ and all $x\in\Omega$. For
convex $\Omega$ and $1<p<(d+2)_+/(d-2)_+$, which includes $p=1+2/d$,
Baras and Cohen \cite{BCo87} proved that bounded nonnegative solutions
blow up completely after $\Tmax$ in this sense; see also
\cite{Ma98,GV97} and \cite[Section~27]{QS19}. By Remark~\ref{rem:proper},
this does not depend on the choice of truncation. If their theorem is
stated for more regular data, fix $t_1\in(0,\Tmax(w_0))$ and
$n\ge\sup_{[0,t_1]}\|w\|_\infty^p$. Then $W_n=w$ on $[0,t_1]$ by
uniqueness, so $W_n(t_1+\cdot)$ is the solution of \eqref{Pn} with the
smooth datum $w(t_1)$, which vanishes on $\partial\Omega$ and satisfies
$\Tmax(w(t_1))=\Tmax(w_0)-t_1$. Hence \eqref{CB} holds if $\Omega$ is
convex.
\end{remark}

\subsection{Proof of Theorem~\ref{cor:OP5}}

\begin{proof}[Proof of Theorem~\ref{cor:OP5}]
Property (i) means that the middle term in \eqref{sandwich} vanishes, and
(iii) that $\Ts(u_0)=0$. By Theorem~\ref{thm:main}(a), (ii) holds if and
only if $\Ts(u_0)=0$: if $\Ts(u_0)>0$, then $U<\infty$ a.e. in
$\Omega\times(0,\Ts(u_0))$. The implications then follow from
\eqref{sandwich}, and the equivalence under \eqref{CB} from
Theorem~\ref{thm:main}(c).

It remains to treat the datum of \cite{BC96}. Brezis and Cazenave fix an
open ball $\omega$ with $\overline\omega\subset\Omega$ and construct
$0\le u_0\in L^1(\Omega)$ with
\begin{equation}\label{BC-claim}
 \int_0^1\!\!\int_\omega\big(S(s)u_0\big)^p\dd x\dd s=+\infty
\end{equation}
\cite[proof of Theorem~11]{BC96}. Suppose that $u$ is an integral
solution on $(0,T)$ for some $T>0$. Choose $(t,x)\in(0,T)\times\Omega$
at which \eqref{integral} holds and $u(t,x)<\infty$, and put
$t':=\min\{t/2,1\}$. Since $t-s\in[t/2,t]$ for $0<s<t'$, continuity and positivity of $G$
give
\[
 c:=\min\big\{G(\sigma,x,y):\ \sigma\in[t/2,t],\ y\in\overline\omega\big\}>0 .
\]
Using \eqref{integral} and $u(s)\ge S(s)u_0$ a.e.,
\[
 \infty>u(t,x)\ge\int_0^{t'}\!\!\int_\omega G(t-s,x,y)u(s,y)^p\dd y\dd s
 \ge c\int_0^{t'}\!\!\int_\omega\big(S(s)u_0\big)^p\dd y\dd s .
\]
For $s\in[t',1]$ we have $\|S(s)u_0\|_\infty\le(4\pi t')^{-d/2}\|u_0\|_1$,
hence
\[
 \int_{t'}^1\!\!\int_\omega\big(S(s)u_0\big)^p\dd x\dd s
 \le|\omega|\,(4\pi t')^{-dp/2}\,\|u_0\|_1^p<\infty .
\]
Together with the previous display, this contradicts
\eqref{BC-claim}. Hence $\Ts(u_0)=0$, and (i), (ii) follow.
\end{proof}

\section{Proof of Theorems~\ref{thm:conc}, \ref{thm:data} and
Corollary~\ref{cor:AR}}\label{sec:conc}

Throughout this section $x_0\in\Omega$, $R_0:=\operatorname{dist}(x_0,\partial\Omega)$
and $R:=R_0/8$, so that $B_{8R}(x_0)\subset\Omega$; the constants are
those of Section~\ref{subsec:constants}.

\subsection{Kaplan's criterion and one step}

We begin with a scale-invariant local version of Kaplan's argument
\cite{Ka63}; for proper minimal solutions a closely related statement is
\cite[Lemma~3.3]{AR04}.

\begin{lemma}[Local Kaplan criterion]\label{lem:kaplan}
Let $0\le w_0\in L^\infty(\Omega)$ with solution $w$, let $x_0\in\Omega$
and $h>0$ with $\overline D\subset\Omega$, where
$D:=B_{b\sqrt h}(x_0)$, and let $\psi>0$ be the first Dirichlet
eigenfunction of $D$, normalised by $\int_D\psi=1$. If
$t_1\in(0,\Tmax(w_0))$ and
\[
 \int_Dw(t_1)\psi\ge M_d\,h^{-d/2},
\]
which holds in particular if $w(t_1)\ge M_dh^{-d/2}$ in $D$, then
$\Tmax(w_0)<t_1+h$.
\end{lemma}

\begin{proof}
We have $-\Delta\psi=(\Lambda/h)\psi$ in $D$. Suppose, for
contradiction, that $\Tmax(w_0)\ge t_1+h$. For $0\le s\le h/2$ put
$Y(s):=\int_Dw(t_1+s)\psi$. Since $\overline D\subset\Omega$ and
$t_1>0$, $w$ is smooth on $[t_1,t_1+h/2]\times\overline D$. Using the
equation and Green's formula,
\[
 Y'(s)=\int_D(\Delta w+w^p)\psi
 =\int_Dw\,\Delta\psi+\int_{\partial D}\big(\psi\,\partial_\nu w
 -w\,\partial_\nu\psi\big)+\int_Dw^p\psi ,
\]
where $\nu$ is the outer normal. The term $\psi\,\partial_\nu w$
vanishes because $\psi=0$ on $\partial D$, and $-w\,\partial_\nu\psi\ge0$
because $w\ge0$ and $\partial_\nu\psi\le0$. Since $\int_D\psi=1$, Jensen's
inequality gives $\int_Dw^p\psi\ge Y^p$. Hence
\[
 Y'(s)\ge-\frac\Lambda hY(s)+Y(s)^p .
\]
We rescale to unit time. The function $Z(\sigma):=h^{d/2}Y(h\sigma)$,
$0\le\sigma\le\frac12$, is $C^1$, $Z(0)\ge M_d$, and
\[
\begin{aligned}
 Z'(\sigma)=h^{\frac d2+1}Y'(h\sigma)
 &\ge-\Lambda\,h^{\frac d2}Y(h\sigma)+h^{\frac d2+1}Y(h\sigma)^p\\
 &=-\Lambda Z(\sigma)+h^{\frac d2+1-\frac{dp}2}Z(\sigma)^p
 =-\Lambda Z(\sigma)+Z(\sigma)^p,
\end{aligned}
\]
since $\frac d2+1-\frac{dp}2=0$ for $p=1+\frac2d$. Let
\[
 \sigma_1:=\sup\big\{\sigma\in[0,\tfrac12]:\ Z\ge Z(0)\text{ on }[0,\sigma]\big\}.
\]
On $[0,\sigma_1]$ we have $Z^{p-1}\ge M_d^{p-1}=M_d^{2/d}\ge2\Lambda$,
hence
\[
 Z'\ge Z^p-\Lambda Z=Z^p\big(1-\Lambda Z^{1-p}\big)\ge\tfrac12Z^p>0 ;
\]
thus $Z$ is increasing on $[0,\sigma_1]$ and, by continuity,
$\sigma_1=\frac12$. Consequently
$(Z^{1-p})'=(1-p)Z^{-p}Z'\le-\frac{p-1}2$ on $[0,\frac12]$, and
integrating over $[0,\frac12]$,
\[
 0<Z(\tfrac12)^{1-p}\le Z(0)^{1-p}-\frac{p-1}4\le M_d^{-2/d}-\frac1{2d}\le0,
\]
since $p-1=2/d$ and $M_d^{2/d}\ge2d$. This contradiction proves the
lemma.
\end{proof}

The part $(2\Lambda)^{d/2}h^{-d/2}=(2\lambda_1(D))^{d/2}$ of the
threshold is the one of \cite[Lemma~3.3]{AR04}. Since
$\psi\le\|\psi_{B_1}\|_\infty(b\sqrt h)^{-d}$, where $\psi_{B_1}$ is the
normalised eigenfunction of $B_1$, the hypothesis of
Lemma~\ref{lem:kaplan} forces $\int_Dw(t_1)\ge M_db^d/\|\psi_{B_1}\|_\infty$,
a mass bounded below by a constant depending only on $d$, whatever the
radius of $D$.

\begin{lemma}[One step]\label{lem:onestep}
Let $0\le w_0\in L^\infty(\Omega)$ with solution $w$, let $\rho>0$,
$\mu\ge0$ with $\int_{B_\rho(x_0)}w_0\ge\mu$, and let
$h\in[16\rho^2,2^{-10}R_0^2]$ with $\Tmax(w_0)>h$. Then
\[
 w(h,x)\ge h^{-d/2}\Big[c_1\,\mu+c_2\,\mu^p\log\frac{h}{16\rho^2}\Big]
 \qquad\text{for all }x\in B_{b\sqrt h}(x_0).
\]
\end{lemma}

\begin{proof}
Note that $2^{-10}R_0^2=R^2/16$. All points below lie in $B_R(x_0)$,
since $\rho\le\sqrt h/4\le R/16$, and all times lie in $(0,R^2]$. For
$x,y\in B_{4R}(x_0)$ the segment joining them stays at distance at least
$4R$ from $\partial\Omega$, so \eqref{vdB} with $\ell=4R$ gives
\begin{equation}\label{kernel-lower}
 G(t,x,y)\ge c_0\,t^{-d/2}\,e^{-|x-y|^2/(4t)}
 \qquad\text{for }x,y\in B_{4R}(x_0),\ 0<t\le R^2.
\end{equation}
In each of the three estimates below the Gaussian exponent in
\eqref{kernel-lower} is at most $\frac14$.

\emph{Linear part.} For $x\in B_{b\sqrt h}(x_0)$ and
$z\in B_\rho(x_0)$ we have $|x-z|\le(b+\frac14)\sqrt h$ and
$(b+\frac14)^2/4\le\frac14$, so that, by \eqref{kernel-lower} and
$c_1=c_0e^{-1/4}$,
\[
 S(h)w_0(x)\ge\int_{B_\rho(x_0)}G(h,x,z)\,w_0(z)\dd z
 \ge c_0h^{-d/2}e^{-1/4}\int_{B_\rho(x_0)}w_0\ge c_1\,\mu\,h^{-d/2}.
\]

\emph{Lower bound at intermediate times.} Let $4\rho^2\le s\le h/4$ and
$y\in B_{a\sqrt s}(x_0)$. For $z\in B_\rho(x_0)$ we have
$\rho\le\sqrt s/2$, hence $|y-z|\le(a+\frac12)\sqrt s$, and
$(a+\frac12)^2/4\le\frac14$. Therefore, as in the linear part,
\begin{equation}\label{Ss-lower}
 S(s)w_0(y)\ge\int_{B_\rho(x_0)}G(s,y,z)\,w_0(z)\dd z
 \ge c_1\,\mu\,s^{-d/2}.
\end{equation}

\emph{Duhamel part.} For $x\in B_{b\sqrt h}(x_0)$,
$y\in B_{a\sqrt s}(x_0)$ and $s\le h/4$, we have
$|x-y|\le b\sqrt h+a\sqrt s\le(b+\frac a2)\sqrt h$ and
$h-s\in[\frac34h,h]$, hence
\[
 \frac{|x-y|^2}{4(h-s)}\le\frac{(b+\frac a2)^2h}{3h}\le\frac14,
 \qquad
 G(h-s,x,y)\ge c_0(h-s)^{-d/2}e^{-1/4}\ge c_1h^{-d/2}.
\]
By
\eqref{mild} and $w(s)\ge S(s)w_0$,
\[
 w(h,x)\ge S(h)w_0(x)+\int_{4\rho^2}^{h/4}\int_{B_{a\sqrt s}(x_0)}
 G(h-s,x,y)\big(S(s)w_0(y)\big)^p\dd y\dd s .
\]
By the bound on $G(h-s,x,y)$, \eqref{Ss-lower} and
$|B_{a\sqrt s}|=\omega_da^ds^{d/2}$, the double integral is at least
\[
\begin{aligned}
 \int_{4\rho^2}^{h/4}c_1h^{-d/2}\,\omega_da^ds^{\frac d2}
 \big(c_1\mu s^{-\frac d2}\big)^p\dd s
 &=c_1^{1+p}\omega_da^d\,\mu^p\,h^{-d/2}
 \int_{4\rho^2}^{h/4}s^{\frac d2-\frac{dp}2}\dd s\\
 &=c_2\,\mu^p\,h^{-d/2}\int_{4\rho^2}^{h/4}\frac{\dd s}s
 =c_2\,\mu^p\,h^{-d/2}\log\frac h{16\rho^2},
\end{aligned}
\]
because $\frac d2-\frac{dp}2=-1$ when $p=1+2/d$. Adding the linear part
gives the claim.
\end{proof}

\subsection{Proof of Theorem~\ref{thm:conc}}

At the critical exponent the time integral in Lemma~\ref{lem:onestep}
counts the parabolic scales between $\rho^2$ and $h$. In the proof below,
the scales $\rho_k$, masses $\mu_k$, waiting times $h_k$ and times $t_k$
evolve as follows:
\begin{center}
\begin{tabular}{@{}cccc@{}}
\toprule
mass & radius & waiting time & time\\
\midrule
$\mu_k=2^k\mu$ & $\rho_{k+1}=b\sqrt{h_k}$ &
$h_k=16\rho_k^2\exp\big(\frac2\beta\,\mu_k^{1-p}\big)$ & $t_{k+1}=t_k+h_k$\\
\bottomrule
\end{tabular}
\end{center}
The mass doubles at each step, and the waiting times shrink geometrically
in the exponent, so the total time is of order $\rho^2\exp(C\mu^{-2/d})$.

\begin{proof}[Proof of Theorem~\ref{thm:conc}]
\emph{Step 1: bounded data.} We first show the following. \emph{Let
$0\le w_0\in L^\infty(\Omega)$, $\rho>0$ and $\mu>0$ with
$\int_{B_\rho(x_0)}w_0\ge\mu$, and suppose that
$\tau:=C_d\rho^2\exp(C_d\mu^{-2/d})<R_0^2$. Then
$\Tmax(w_0)<2^{-10}\tau$.}

Let $K:=0$ if $\mu\ge\mu_1$, and $K:=\lceil\log_2(\mu_1/\mu)\rceil$
otherwise. Define, for $0\le k<K$,
\[
 \mu_k:=2^k\mu,\qquad L_k:=\frac2\beta\,\mu_k^{1-p},\qquad
 h_k:=16\rho_k^2e^{L_k},\qquad \rho_{k+1}:=b\sqrt{h_k},\qquad
 t_{k+1}:=t_k+h_k,
\]
with $\rho_0:=\rho$ and $t_0:=0$, and put $\mu_K:=2^K\mu\ge\mu_1$,
$h_K:=16\rho_K^2$ and $t_{K+1}:=t_K+2h_K$.

\emph{Step 1a: bookkeeping.} Let $S:=\sum_{j<K}L_j$. Since $1-p=-2/d$,
\[
 S\le\frac2\beta\,\mu^{-2/d}\sum_{j\ge0}2^{-2j/d}=C_1\mu^{-2/d}.
\]
Since $16b^2\le1$, we have
\[
 \rho_{k+1}^2=b^2h_k=16b^2\rho_k^2e^{L_k}\le\rho_k^2e^{L_k},
\]
and by induction on $k$,
\[
 \rho_k^2\le\rho^2e^{L_0+\dots+L_{k-1}}\le\rho^2e^S\qquad(0\le k\le K).
\]
Hence
\[
 h_k=16\rho_k^2e^{L_k}\le16\rho^2e^{L_0+\dots+L_k}\le16\rho^2e^S\quad(k<K),
 \qquad
 h_K=16\rho_K^2\le16\rho^2e^S .
\]
Next we bound $K$. If $\mu\ge\mu_1$, then $K=0$. If $\mu<\mu_1$, then
$K\le1+\log_2(\mu_1/\mu)$; since $\log y\le y$ for $y>0$, applied with
$y=x^{2/d}$,
\[
 \log_2x=\frac d{2\log2}\,\log x^{2/d}\le\frac d{2\log2}\,x^{2/d}
 \qquad(x\ge1),
\]
and with $x=\mu_1/\mu$ this gives $\log_2(\mu_1/\mu)\le C_2\mu^{-2/d}$.
In both cases, using $3+y\le3e^y$ for $y\ge0$,
\[
 K+2\le3+C_2\mu^{-2/d}\le3e^{C_2\mu^{-2/d}} .
\]
Therefore, since $C_1+C_2\le C_d$ and $48\le2^{-10}C_d$ by \eqref{CN},
\[
\begin{aligned}
 t_{K+1}&=\sum_{k<K}h_k+2h_K\le16(K+2)\rho^2e^S
 \le48\rho^2e^{(C_1+C_2)\mu^{-2/d}}\\
 &\le2^{-10}C_d\rho^2e^{C_d\mu^{-2/d}}=2^{-10}\tau<2^{-10}R_0^2 .
\end{aligned}
\]
In particular $16\rho_k^2\le h_k\le2^{-10}R_0^2$ for every $k\le K$.

\emph{Step 1b: iteration.} We claim that, for $0\le k\le K$, either
$\Tmax(w_0)\le t_k$, or $t_k<\Tmax(w_0)$ and
$\int_{B_{\rho_k}(x_0)}w(t_k)\ge\mu_k$. This holds for $k=0$. Let
$k<K$ and assume $t_k<\Tmax(w_0)$ and
$\int_{B_{\rho_k}(x_0)}w(t_k)\ge\mu_k$. If $\Tmax(w_0)\le t_{k+1}$ there
is nothing to prove. Otherwise, by Lemma~\ref{lem:bounded}(e), we may
apply Lemma~\ref{lem:onestep} to the datum $w(t_k)$ with $\rho_k$,
$\mu_k$ and $h=h_k$. Since $\log(h_k/16\rho_k^2)=L_k$, keeping only the Duhamel part we get
\[
 w(t_{k+1},x)\ge c_2\,\mu_k^p\,L_k\,h_k^{-d/2}
 \qquad\text{for }x\in B_{b\sqrt{h_k}}(x_0)=B_{\rho_{k+1}}(x_0).
\]
Integrating over this ball, whose volume is $\omega_db^dh_k^{d/2}$, we
obtain
\[
 \int_{B_{\rho_{k+1}}(x_0)}w(t_{k+1})\ge\omega_db^dc_2\,\mu_k^pL_k
 =\beta\mu_k^p\cdot\frac2\beta\mu_k^{1-p}=2\mu_k=\mu_{k+1}.
\]

\emph{Step 1c: conclusion.} If $\Tmax(w_0)\le t_K+h_K$, then
$\Tmax(w_0)<t_{K+1}\le2^{-10}\tau$. Otherwise, by the claim and
Lemma~\ref{lem:onestep}, applied to $w(t_K)$ with $h=h_K=16\rho_K^2$ and
the linear part only,
\[
 w(t_K+h_K,x)\ge h_K^{-d/2}c_1\mu_K\ge M_dh_K^{-d/2}
 \qquad\text{for }x\in B_{b\sqrt{h_K}}(x_0),
\]
since $c_1\mu_K\ge c_1\mu_1=M_d$. As $b\sqrt{h_K}<R_0$,
Lemma~\ref{lem:kaplan} with $t_1:=t_K+h_K>0$ gives
$\Tmax(w_0)<t_K+2h_K=t_{K+1}\le2^{-10}\tau$.

\emph{Step 2: general data.} We pass from bounded data to $u_0$ with the
concave transformation. Lemma~\ref{lem:concave} requires a datum
$w\le\lambda u_0$ with $\lambda<1$, and the truncation $\min\{u_0,n\}$
loses part of the mass in $B_\rho(x_0)$; each of these two effects costs
a factor $\theta^{1/2}$, which is why $\mu$ is replaced by $\theta\mu$
below. Let $u_0$ be as in the theorem and suppose,
for contradiction, that $\Ts(u_0)>\tau$. Since
$\tau<\min\{\Ts(u_0),R_0^2\}$, we can fix $\theta\in(0,1)$ such that
\[
 \tau_\theta:=C_d\rho^2\exp\big(C_d(\theta\mu)^{-2/d}\big)
 <\min\{\Ts(u_0),R_0^2\}.
\]
By monotone convergence there is $n$ with
$\int_{B_\rho(x_0)}\min\{u_0,n\}\ge\theta^{1/2}\mu$. The bounded datum
$w:=\theta^{1/2}\min\{u_0,n\}\le\theta^{1/2}u_0$ satisfies
$\int_{B_\rho(x_0)}w\ge\theta\mu$, so Step~1 gives $\Tmax(w)<\tau_\theta<\Ts(u_0)$. This contradicts
Lemma~\ref{lem:concave}. Hence $\Ts(u_0)\le\tau$. The assertions on
$\Tmax$ and on $U$ follow from Remark~\ref{rem:classes}(1) and
Theorem~\ref{thm:main}(a).

\emph{Step 3: the equivalent form.} For $C_d\rho^2<T'$, raising both sides of \eqref{necessary} to the
power $-2/d$ shows that \eqref{necessary} is equivalent to
\[
 C_d\,\mu^{-2/d}\ge\log\frac{T'}{C_d\rho^2},
 \qquad\text{that is,}\qquad
 \tau=C_d\rho^2\exp\big(C_d\mu^{-2/d}\big)\ge T' .
\]
If \eqref{P} has an integral solution on $(0,T)$, then $\Ts(u_0)\ge T$,
so $\tau<T'=\min\{T,R_0^2\}$ is impossible by Steps~1--2; hence
\eqref{necessary} holds. Conversely, if \eqref{necessary} holds for every
$T<\Ts(u_0)$ and $\tau<R_0^2$, then $\tau\ge\min\{T,R_0^2\}$ for every such
$T$, so $\Ts(u_0)\le\tau$.
\end{proof}

\begin{remark}[One step]\label{rem:onestep}
A single application of Lemma~\ref{lem:onestep} with
\[
 \log\frac h{16\rho^2}=\frac{M_d}{c_2\mu^p},
 \qquad\text{that is,}\qquad
 h=16\rho^2\exp\Big(\frac{M_d}{c_2\mu^p}\Big),
\]
gives $w(h)\ge M_dh^{-d/2}$ in $B_{b\sqrt h}(x_0)$, and
Lemma~\ref{lem:kaplan} with $t_1=h$ then gives
\[
 \Tmax(w_0)<2h=32\rho^2\exp\Big(\frac{M_d}{c_2\mu^p}\Big)
\]
whenever the right-hand side is at most $2^{-9}R_0^2$. This weaker estimate has the exponent $-p$
instead of $-2/d$. It already suffices for the qualitative part of
Theorem~\ref{thm:data}(a) under the stronger hypothesis
$(\log\frac1r)^{1/p}m(r)\to\infty$, for instance for $u_{\gamma,A}$ with
$\gamma<1+\frac d{d+2}$. The iteration over scales is what produces the
optimal exponent.
\end{remark}

\begin{remark}[Optimality]\label{rem:optimal}
Let $0\le\psi_0\in C_c^\infty(B_1)$ with $\int\psi_0=1$, let
$B_\rho(x_0)\subset\Omega$, and let
$w_0:=\mu\rho^{-d}\psi_0((\cdot-x_0)/\rho)$. Since
$G_\Omega\le G_{\R^d}$, the solution of the Cauchy problem with datum
$w_0$ (extended by zero) is a supersolution of \eqref{P}, so $\Tmax(w_0)$
is at least the lifespan of the Cauchy problem. By the $L^1$-invariant
scaling the latter equals $\rho^2$ times the lifespan for the datum
$\mu\psi_0$, which for small $\mu$ is at least $\exp(c\mu^{-2/d})$ by the
critical lifespan estimates of Lee and Ni \cite{LN92}. Hence
$\Ts(w_0)\ge\Tmax(w_0)\ge\rho^2\exp(c\mu^{-2/d})$ for small $\mu$, and
the exponent in Theorem~\ref{thm:conc} cannot be improved.
\end{remark}

\subsection{Proof of Theorem~\ref{thm:data}}\label{sec:data}

\begin{proof}[Proof of Theorem~\ref{thm:data}\textnormal{(a)}]
Fix $K$ with
\[
 \mu_d<K<\limsup_{r\downarrow0}\Big(\log\frac1r\Big)^{d/2}m(r),
\]
and choose $r_j\downarrow0$ with $r_j<1$ and
\[
 m(r_j)\ge K\Big(\log\frac1{r_j}\Big)^{-d/2}\qquad\text{for all }j.
\]
Fix $\theta\in(\mu_d/K,1)$ and put
\[
\begin{aligned}
 \tau_j&:=C_d r_j^2\exp\Big(C_d\Big(\theta K
 \Big(\log\frac1{r_j}\Big)^{-d/2}\Big)^{-2/d}\Big)\\
 &=C_d r_j^2\exp\Big(C_d(\theta K)^{-2/d}\log\frac1{r_j}\Big)
 =C_d\,r_j^{\,2-C_d(\theta K)^{-2/d}} .
\end{aligned}
\]
Since $\theta K>\mu_d=(C_d/2)^{d/2}$, we have
$C_d(\theta K)^{-2/d}<C_d\mu_d^{-2/d}=2$, so $\tau_j\to0$. Let $w_n\to u_0$ in $L^1(B_\delta(x_0))$. For $j$ with
$r_j<\delta$ and $\tau_j<R_0^2$ there is $n_j$ such that
\[
 \int_{B_{r_j}(x_0)}w_n\ge\theta\,m(r_j)
 \ge\theta K\Big(\log\frac1{r_j}\Big)^{-d/2}\qquad(n\ge n_j),
\]
and Theorem~\ref{thm:conc} with $\rho=r_j$ gives $\Ts(w_n)\le\tau_j$ for
$n\ge n_j$. Hence $\limsup_n\Ts(w_n)\le\tau_j$ for all large $j$, that
is, $\Ts(w_n)\to0$. Taking $w_n=u_0$ gives $\Ts(u_0)=0$, and (i), (ii)
follow from Theorem~\ref{cor:OP5}.
\end{proof}

\begin{proof}[Proof of Theorem~\ref{thm:data}\textnormal{(b)}]
Write $u:=u_{\gamma,A}$. The choice of $\delta$ ensures
$\log\frac1r>1+\gamma/d$ for $0<r<\delta$, so $u$ is radially decreasing
on $B_\delta(x_0)\setminus\{x_0\}$. For $0<r\le\delta$,
\[
 m(r)=A\sigma_d\int_0^r\frac{\dd s}{s\,(\log\frac1s)^\gamma}
 =\frac{A\sigma_d}{\gamma-1}\Big(\log\frac1r\Big)^{1-\gamma},
\]
so $u\in L^1$, and
\[
 \Big(\log\frac1r\Big)^{d/2}m(r)
 =\frac{A\sigma_d}{\gamma-1}\Big(\log\frac1r\Big)^{\frac d2+1-\gamma}.
\]
This tends to $\infty$ if $\gamma<1+d/2$, and equals $2A\sigma_d/d$ if
$\gamma=1+d/2$.

\emph{The case $1<\gamma<1+d/2$.} By part (a), $u$ satisfies (i)--(iii).
The lower bound in \eqref{rate} is Lemma~\ref{lem:bounded}(f). For the
upper bound, let $n$ be large and let $r_n\in(0,\delta)$ satisfy $u=n$ on
$\partial B_{r_n}(x_0)$, so that $\min\{u,n\}=u$ on
$\{r_n\le|x-x_0|<\delta\}$. Put $L:=\log\frac1{r_n}\ge1$. The condition
$u=n$ on $\partial B_{r_n}(x_0)$ reads $Ar_n^{-d}L^{-\gamma}=n$; taking
logarithms,
\[
 \log A+dL-\gamma\log L=\log n,
 \qquad\text{that is,}\qquad
 L=\frac1d\big(\log n-\log A+\gamma\log L\big),
\]
hence
\[
 -2L=-\frac2d\log n+\frac2d\log A-\frac{2\gamma}d\log L
 \le-\frac2d\log n+\frac2d\log A ,\qquad L\le\frac2d\log n
 \quad\text{for large }n.
\]
Let $\eps\in(0,\frac12]$,
$r:=r_n^{1-\eps}$ and $\ell:=\log\frac1r=(1-\eps)L$. Then
\[
 \int_{B_r(x_0)}\min\{u,n\}\ge\int_{r_n<|x-x_0|<r}u
 =\frac{A\sigma_d}{\gamma-1}\big(1-(1-\eps)^{\gamma-1}\big)\ell^{1-\gamma}
 \ge c_*\,\eps\,\ell^{1-\gamma},
\]
with $c_*:=A\sigma_d\min\{1,(\gamma-1)^{-1}\}$, since
$1-(1-\eps)^{a}\ge\min\{1,a\}\eps$ for $a>0$. Theorem~\ref{thm:conc} with
$\rho=r$ and this mass gives, as long as the resulting $\tau$ is smaller
than $R_0^2$,
\[
 \log\Tmax(\min\{u,n\})\le\log C_d-2(1-\eps)L
 +C_d(c_*\eps)^{-2/d}L^{\vartheta},
 \qquad\vartheta:=\frac{2(\gamma-1)}d\in(0,1).
\]
Writing $-2(1-\eps)L=-2L+2\eps L$, the price of the choice of $r$ is the
term $2\eps L$, and the gain is that the last term contains
$\eps^{-2/d}L^\vartheta$ instead of $L$. We balance the two by requiring
$\eps L=\eps^{-2/d}L^\vartheta$, that is,
\[
 \eps:=L^{-\frac{(1-\vartheta)d}{d+2}}\le\tfrac12\quad\text{(for large $n$)},
 \qquad\text{so that}\qquad
 \eps L=\eps^{-2/d}L^\vartheta=L^{\varsigma},
\]
where
\[
 \varsigma:=1-\frac{(1-\vartheta)d}{d+2}=\frac{2+\vartheta d}{d+2}
 =\frac{2\gamma}{d+2}\in(0,1).
\]
With this choice
$2\eps L+C_d(c_*\eps)^{-2/d}L^\vartheta=\big(2+C_dc_*^{-2/d}\big)L^\varsigma$,
so the right-hand side above equals
\[
 \log C_d-2L+\big(2+C_dc_*^{-2/d}\big)L^{\varsigma}\longrightarrow-\infty
 \qquad(n\to\infty),
\]
since $\varsigma<1$ and $L\to\infty$. Hence the condition $\tau<R_0^2$
holds for large $n$, and by the bounds on $L$ above
\[
\begin{aligned}
 \log\Tmax(\min\{u,n\})
 &\le\log C_d+\frac2d\log A-\frac2d\log n
 +\big(2+C_dc_*^{-2/d}\big)\Big(\frac2d\log n\Big)^{\varsigma}\\
 &\le-\frac2d\log n+C''(\log n)^{\frac{2\gamma}{d+2}}
 \qquad(n\ge n_0),
\end{aligned}
\]
which is the upper bound in \eqref{rate}.

\emph{The case $\gamma=1+d/2$.} Put $A_c:=\kappa_*(u_{\gamma,1})$, so that
$u_{\gamma,A}=Au_{\gamma,1}$ satisfies (i)--(iii) for $A>A_c$ and (i)
fails for $0<A<A_c$ by Theorem~\ref{thm:main}(d). Since $u_{\gamma,A}$
satisfies (iii) for $A>A_d$ by part (a), $A_c\le A_d$. Finally,
$A_c>0$ because, by \cite{HI18,HI24}, \eqref{P} is locally solvable for
$u_{\gamma,A}$ with $A$ small; for the Cauchy problem this transfers to
$\Omega$ as in the proof of part (c).
\end{proof}

\begin{proof}[Proof of Theorem~\ref{thm:data}\textnormal{(c)}]
Let $q\in[0,d/2)$ and choose $\gamma\in(1+q,1+d/2)$. On $B_\delta(x_0)$
we have $u_{\gamma,1}\le|x-x_0|^{-d}$ and
$\log(e+|x-x_0|^{-d})\le(d+1)\log\frac1{|x-x_0|}$, hence
\[
 u_{\gamma,1}[\log(e+u_{\gamma,1})]^q\le(d+1)^q|x-x_0|^{-d}
 \Big(\log\frac1{|x-x_0|}\Big)^{q-\gamma},
\]
which is integrable because $\gamma-q>1$. By (b), $u_{\gamma,1}$
satisfies (i)--(iii).

Conversely, let $u_0\ge0$ with
$\int_\Omega u_0[\log(e+u_0)]^{d/2}<\infty$. Since
$\log(e+2s)\le2\log(e+s)$ for $s\ge0$, also
$\int_\Omega2u_0[\log(e+2u_0)]^{d/2}<\infty$. By \cite[Theorem~1.3]{Mi21},
the Cauchy problem whose datum is the zero extension of $2u_0$ has a
nonnegative solution $\bar v$ on $(0,T)$ for some $T>0$.
Since $G_\Omega\le G_{\R^d}$, the restriction of $\bar v$ to $\Omega$
satisfies \eqref{super} with datum $2u_0$, so by Lemma~\ref{lem:U} the
minimal solution with datum $2u_0$ is finite a.e. in $\Omega\times(0,T)$,
hence an integral solution there. Thus $\Ts(u_0)\ge\Ts(2u_0)\ge T>0$, so
(iii) and (ii) fail, and (i) fails by \eqref{sandwich} with $\kappa=2$.
Sufficient conditions directly on domains are given in
\cite[Section~5]{HI24}.
\end{proof}

\subsection{Proof of Corollary~\ref{cor:AR}}

\begin{proof}[Proof of Corollary~\ref{cor:AR}]
Let $0\le v_n\in C_c^\infty(\Omega)$ converge to $\nu$ narrowly, and in
$L^1(\Omega)$ if $\nu\in L^1(\Omega)$; such $v_n$ are obtained by
restricting $\nu$ to compact sets exhausting $\Omega$ and mollifying.
Fix $x_0\in\Omega$, let $R_0:=\operatorname{dist}(x_0,\partial\Omega)$, and
let $\psi_0$ be as in Remark~\ref{rem:optimal}. Choose $\rho_n\in(0,R_0)$
so small that
\[
 C_d\,\rho_n^2\exp\big(C_dn^{2/d}\big)<\min\Big\{\frac1n,\,R_0^2\Big\},
\]
and put
\[
 w_n:=v_n+\frac1n\,\rho_n^{-d}\psi_0\Big(\frac{\cdot-x_0}{\rho_n}\Big)
 \in C_c^\infty(\Omega).
\]
The added bumps are nonnegative and have mass $1/n\to0$, so $w_n\to\nu$
narrowly, and in $L^1$ if $\nu\in L^1$. Since
$\int_{B_{\rho_n}(x_0)}w_n\ge1/n$, Theorem~\ref{thm:conc} with $\rho=\rho_n$
and $\mu\ge1/n$ gives $\Ts(w_n)\le C_d\rho_n^2\exp(C_dn^{2/d})<1/n$. The
remaining assertions follow from Remark~\ref{rem:classes}(1) and
Theorem~\ref{thm:main}(a).
\end{proof}

\begin{remark}[Nested spikes]\label{rem:spikes}
A self-contained example is $u_0:=\sum_jv_j$ with
$v_j:=2^{-j}r_j^{-d}\psi_0((\cdot-x_0)/r_j)$, $\psi_0$ as in
Remark~\ref{rem:optimal}, and $r_j$ so small that
\[
 C_d\,r_j^2\exp\big(C_d\,2^{2j/d}\big)<\min\{2^{-j},R_0^2\}.
\] Since
$\int_{B_{r_j}(x_0)}u_0\ge2^{-j}$, Theorem~\ref{thm:conc} gives
$\Ts(u_0)<2^{-j}$ for every $j$; hence $u_0$ satisfies (i)--(iii). This
is the construction of Corollary~\ref{cor:AR} with summable instead of
vanishing masses. For statements on the Orlicz scale, the explicit data
of Theorem~\ref{thm:data}(b) are more convenient.
\end{remark}


\section*{Acknowledgments}
The author was supported by the Science Committee of the Ministry of Science and Higher Education of the Republic of Kazakhstan (Grant No. BR31714735).

\section*{Declaration on the use of generative AI}

The author used Claude (Anthropic) for language editing, improving the presentation of the manuscript, and assisting in the verification of mathematical arguments. All AI-assisted material was carefully reviewed and independently verified by the author. The author takes full responsibility for the content, correctness, and integrity of the manuscript.

\section*{Declaration of competing interest}
The author declare that there is no conflict of interest.

\section*{Data Availability Statements}
The manuscript has no associated data.

\end{document}